\documentclass [12pt,a4paper,reqno]{amsart}
\usepackage{amsmath,amsthm}
\usepackage{amscd}
\usepackage{amsfonts}
\input amssymb.sty
\DeclareMathOperator{\Hom}{Hom}

\newcommand{\ef}{\end{equation}}
\chardef\bslash=`\\ 

\newtheorem{thm}{Theorem}
\newtheorem*{thm*}{Theorem}
\newtheorem*{conjecture*}{Conjecture}

\theoremstyle{definition}

\newtheorem*{remark*}{Remarks}
\newtheorem*{examples*}{Examples}
\newtheorem*{defn*}{Definition}
\newtheorem*{cor*}{Corollary}

\newtheorem{rem}{Remark}

\numberwithin{equation}{section}

\newcommand{\Th}{\Theta}
\newcommand{\C}{\mathcal{C}}
\newcommand{\V}{\mathcal{V}}
\newcommand{\Q}{\mathbf{Q}}
\newcommand{\st}{\sigma}

 \renewcommand{\sectionmark}[1]{}

 \usepackage{amsfonts}

\newcommand{\vp}{\varphi}

\renewcommand{\a}{\alpha}

 \date{}
\begin{document}

\title{A generalization and a new proof of Plotkin's reduction theorem}
\author[Grigori Zhitomirski]{G. Zhitomirski\\
 Department of Mathematics and Statistics \\
 Bar-Ilan University, 52900 Ramat Can, Israel}
\thanks {This research is partially supported by THE ISRAEL SCIENCE FOUNDATION
 founded by The Israel Academy of Sciences and Humanities - Center of
 Excellence Program.}

\maketitle

\begin{abstract}
It is known that Plotkin's reduction theorem is very important for his theory of universal algebraic geometry
\cite{1,2}. It turns out that this theorem can be generalized to arbitrary categories containing two special
objects and in this case its proof  becomes considerable more simple. This new proof and applications are the
subject of the present paper.
\end{abstract}

\baselineskip 20pt
\bigskip

\centerline{INTRODUCTION}\label{Intro}

An automorphism $\vp$ of a category $\C$ is called inner if it is isomorphic to the identity functor $Id \sb
{\C}$   in the category of all endofunctors of $\C$.  If an automorphism $\vp$ is inner the object $\vp (A)$ is
isomorphic to $A$ for every $\C$-object $A$. It is known \cite{2} that every automorphism satisfying the last
condition is a composition of two automorphisms the first of which preserves the objects and second one is an
inner automorphism. Thus it is sufficient to consider only such automorphisms $\vp$ that preserve the objects,
i.e. $\vp (A)=A$ for all objects $A$ of $\C $. Further every automorphism satisfying the last condition induces
a permutation on the set $\Hom (A,B)$ for every pair of objects $A,B$, particularly it induces an automorphism
of the monoid $End (A)=\Hom (A,A)$ for every object $A$. Below I cite the Plotkin's Reduction Theorem. The
purpose of this theorem is to reduce the verification of a given automorphism of $\C$ to be inner to the same
problem for the full subcategory of $\C$ determined by two special objects.

Let $\V$ be a variety of universal algebras. Consider the category $\Th $  whose objects are all algebras from
$\ V$ and whose morphisms are all homomorphisms of them. Fix an infinite set $X$. Let $\Th \sp 0 $ be the full
subcategory of $\Th $ defined by all free algebras from $\V$ over finite subsets of the set $X$.

The following conditions are assumed:
\begin{enumerate}
\item[1P)] every object of $\Th \sp 0 $ is a hopfian algebra;

\item[2P)] there exists an object $F\sp 0=F(X\sp 0 )$ in $\Th \sp 0 $ generating the whole variety $\V$.
\end{enumerate}

Fix the object  $F\sb 0 =F(X\sb 0 )$  in $\Th \sp 0 $ generated by a singular set $X\sb 0 = \{x\sb 0 \}$ and the
homomorphism $\nu \sb 0 :F\sp 0 \to F\sb 0 $ induced by the constant map $X\sp 0 \to X\sb 0$ that is $(\forall
x\in X\sp 0 )\; \nu \sb 0 (x)=x\sb 0$.

\begin{thm}\label{theorem1} Let $\vp : \Th \sp 0 \to \Th \sp 0$ be an automorphism which does not change objects. If
$\vp $ induces the identity automorphism on the semigroup  $END (F\sp 0)$ and $\vp (\nu \sb 0 )=\nu \sb 0$, then
$\vp $ is an inner automorphism.
\end{thm}

It should be mention that this theorem was first proved by Berzins [3] for the variety of commutative
associative algebras over an infinite field.

The purpose of the present paper is to show that some hypotheses used in this theorem are unnecessary, the
statement can be generalized and the proof becomes essentially more easy. The author thanks Prof. B. Plotkin for
useful discussions.

\section{A generalized Reduction theorem}
 Let $\C$ be a category with the following conditions:

\begin{enumerate}
\item[1)] there is a faithful functor $\Q: \C \to Set$  such that it is represented by an object $A_0$ ;

\item[2)] there is an object $A^0$ such that for every two objects $A$ and $B$ and for every bijection $s: \Q (A
)\to \Q (B )$ one of the following two dual conditions is satisfied:

if for every morphism $\nu: B \to A\sp 0  $ ($\nu: A\sp 0 \to A $) there exists a morphism $\mu: A\to A \sp 0 $
($\mu: A\sp 0 \to A$) such that
$$ \Q (\mu )= \Q (\nu )\circ s $$
$$(resp. \; \Q (\mu )= s\circ \Q (\nu ))$$
then there exists an isomorphism $\gamma : A \to B$ such that $\Q (\gamma )=s$.
\end{enumerate}

The sense of the last condition using algebraic language is: if the composition of a given bijection of $A$ and
every homomorphism from $A$  to $A\sp 0$ is a homomorphism then this bijection is an isomorphism of $A$.

\begin{thm}\label{theorem2} If $\vp : \C \to \C $ is an automorphism of the category $\C $ that does not change the
objects $A\sb 0 $ and $ A\sp 0 $ and induces the identity map on $Hom (A\sb 0 , A\sp 0) $ then $\vp $ is an
inner automorphism.
\end{thm}
\begin{proof}

Let objects $A\sp 0$ and $A\sb 0$ existing under hypotheses be fixed. And let $u : \Q \to \Hom (A\sb 0 ,?)$ be
an isomorphisms of functors, i.e. a representation of the functor $\Q$. Let $\vp$ be an isomorphism of the
category $\C$ satisfying required conditions. Thus we have for every object $A$ the bijection $\vp \sb A :\Hom
(A\sb 0 ,A)\to \Hom (A\sb 0 ,\vp (A))$ and therefore a bijection $s \sb A $ of the set $\Q (A)$ onto the set $\Q
(\vp(A))$ unique defined by means of the following commutative diagram:

$$
\CD
\Q (A) @> u_A>> \Hom (A\sb 0 , A)\\
@Vs \sb A VV        @V\vp \sb A VV \\
\Q (\vp (A)) @> u_{\vp (A)}>> \Hom (A\sb 0 , \vp (A))
\endCD
$$

It means that $s \sb A = u\sb {\vp (A)} \sp {-1} \circ \vp \sb A \circ u\sb A$. Under hypothesis $\vp \sb {A\sp
0}$ is the identity map, hence we have $s \sb {A \sp 0}=1\sb {\Q (A\sp 0)}$. Further, for every morphism $\nu
:A\to B$ we have two commutative diagrams:

$$
\CD
\Q (A) @> u_A>> \Hom (A\sb 0 , A)\\
@V\Q (\nu ) VV   @V\Hom (A\sb 0 , \nu) VV \\
\Q (B) @> u_B>> \Hom (A\sb 0 , B)
\endCD
$$
and

$$
\CD
\Q (\vp (A)) @> u_{\vp (A)}>> \Hom (A\sb 0 , \vp (A))\\
@V\Q (\vp (\nu )) VV   @V\Hom (A\sb 0 , \vp (\nu)) VV \\
\Q (\vp (B)) @> u_{\vp (B)}>> \Hom (A\sb 0 , \vp (B))
\endCD
$$

Let $\mu : A\sb 0 \to A$. Since $\vp$ is an automorphism we have:
$$ \Hom (A\sb 0 , \vp (\nu))(\mu)=\vp (\nu))\circ \mu =\vp (\nu \circ \vp \sp {-1}(\mu))=
\vp (\Hom (A\sb 0 , \nu)(\vp \sp {-1}(\mu))).$$

Hence $\Hom (A\sb 0 , \vp (\nu))= \vp \circ \Hom (A\sb 0 , \nu)\circ \vp \sp {-1}$ and from those two commutative
diagrams follows the fact:

$$\Q (\vp (\nu ))= u\sb {\vp (B)} \sp {-1} \circ \Hom (A\sb 0 , \vp (\nu)) \circ u\sb {\vp (A)} = $$
$$=u\sb {\vp (B)} \sp {-1} \circ \vp \circ \Hom (A\sb 0 , \nu)\circ \vp \sp {-1} \circ u\sb {\vp (A)} = $$
$$=u\sb {\vp (B)} \sp {-1} \circ \vp \circ u \sb B \circ u \sb B \sp {-1} \circ \Hom (A\sb 0 , \nu)\circ u \sb A
\circ u \sb A \sp {-1} \circ \vp \sp {-1} \circ u\sb {\vp (A)} =$$
$$= s \sb B \circ \Q (\nu )\circ s \sb A \sp {-1}.$$

Thus for every morphism $\nu :A\to B$ we obtain the result:
$$\Q (\vp (\nu ))=s \sb B \circ \Q (\nu )\circ s \sb A \sp {-1}.$$

Particularly, for every morphism $\nu :A\to A\sp 0 $ we have $\Q (\vp (\nu ))=s \sb {A\sp 0} \circ \Q (\nu
)\circ s \sb A \sp {-1}= \Q (\nu ) \circ s \sb A \sp {-1}$. Under hypotheses for the category $\C$, there exists
an isomorphism $\st \sb A :A\to \vp (A)$ such that $s \sb A =\Q (\st \sb A )$. The dual case gives the same
conclusion. And finally we obtain that for every $\nu :A\to B $:

$$\vp (\nu )=\st \sb B \circ \nu \circ \st \sb A \sp {-1}.$$
that ends the proof.
\end{proof}

\begin{thm}Theorem 2 implies Theorem 1.
\end{thm}
\begin{proof}
Consider the category  $\C = \Th \sp 0$ described in Introduction, that is  $\V$ is a variety of universal
algebras,  $\Th $ is the corresponding category and $\Th \sp 0$ is the full subcategory of $\Th $ defined by all
free algebras from $\V$ over finite subsets of an infinite fixed set $X$. Let $\Q :\C \to Set $ be the forgetful
functor. It is faithful. Let $F\sb 0$ be a free algebra in $\Th \sp 0$ over an one-element set. Then the functor
$Q$ is represented by the object $F\sb 0$ and the condition (1) for the category $\C$ is satisfied. We write
simply  $A$ instead of $\Q (A)$.

Let $\Th \sp 0$ contain a finitely generated free  algebra $ F \sp 0$ generating the variety $\V$ (condition
(2P)). We prove that the condition (2) for the category $\C$ is also satisfied.  Consider two objects $A$ and
$F$ of $\Th \sp 0$. Let $F=F(x\sb 1, ...,x\sb n)$ and $s: A\to F$ be a bijection with the following condition:
for every homomorphism $\nu: F \to F\sp 0  $ the composition $ \mu = \nu \circ s$ is also homomorphism. Let $o$
be a k-ary operation. For every $k$ elements $a_1,...,a_k \in A$, consider two words in $F$: $u=s(o(a_1,...,a_k
))$ and $v=o(s(a_1 ),...,s(a_k ))$. We have $\nu (u)=\mu(o (a_1,...,a_k))=(o (\mu (a_1),...,\mu (a_k))=\nu (o
(s(a_1 ),...,s(a_k))=\nu (v)$. Since all identities which are satisfied in $ F \sp 0$ are satisfied in $F$ too,
$u=v$. It means that $s$ is also a homomorphism. Thus the condition (2) for $\C$ is satisfied.

Now let an automorphism $\vp : \Th \sp 0 \to \Th \sp 0$ satisfy the hypothesis of Theorem~1. Let $\mu : F\sb 0
\to F\sp 0$ be an arbitrary morphism and let $\nu = \mu \circ \nu \sb 0 : F\sp 0 \to F \sp 0$. Under hypotheses
of Theorem~1 we have:
$$ \nu = \vp (\nu )= \vp (\mu )\circ \vp (\nu \sb 0 )=\vp (\mu ) \circ \nu \sb 0.$$
Hence $\mu \circ \nu \sb 0 = \vp (\mu ) \circ \nu \sb 0$ and therefore $\vp (\mu )= \mu $. Thus the automorphism
$\vp : \Th \sp 0 \to \Th \sp 0$  satisfies the hypothesis of Theorem~2.  Since the statements of Theorem~2 and
Theorem~1 coincide  we deduce Theorem~1 from Theorem~2.

\end{proof}

\begin{rem}
Let $A$ and $B$ be two different objects of a category $\C$ and let $\vp$ be an automorphism of $\C$ satisfying
the condition that there are two isomorphisms $\a : A \to \vp (A)$  and $\beta :B \to \vp (B)$ such that for
every $\nu :A\to B $ we have $\vp (\nu ) = \beta \circ \nu \circ \a \sp {-1}$. Then $\vp $ is a composition of
two automorphisms $\vp = \psi \circ \gamma $ where $\gamma $ is an inner automorphism and $\psi $ preserves
objects $A,B$ and acts as identity on $\Hom (A,B)$. Thus for the categories satisfying conditions (1) and (2) we
obtain from Theorem~2 a necessary and sufficient condition for an automorphism to be inner.
\end{rem}

\begin{rem} The condition (P1) has not been used in the proof given above. Thus it is not necessary for the
Theorem 1. However this condition is used in \cite{2} to prove the fact that every automorphism of the category
$\Th \sp 0 $ takes every object to an isomorphic one. But we have seen that it also is not necessary.

\end{rem}

Let $\C$ be a category of universal algebras and their homomorphisms. To apply Theorem~2 it is necessary to find
two objects $A\sb 0 $ and $A\sp 0 $ satisfying conditions (1) and (2). A functor $\Q $ in such a case is usually
the forgetful functor. If a monogenic free algebra is an object of $\C$ the condition (1) is realized. The
condition 2 is realized if our category contains an object $A\sp 0$  satisfying at least one of the following
two conditions:

a) for every $\C$-object $A$ and every two its different elements $a\sb 1$ and $a\sb 2$ there exists a
homomorphism $\nu : A\to A\sp 0$ with $\nu (a\sb 1) \not = \nu (a\sb 2)$

or

b) for every $\C$-object $A$ and every finite subset $X$ of $A$ having as many elements as  arity of an
operation, there exists a homomorphism $\nu : A\sp 0\to A$ with $X\subset \nu (A\sp 0)$ .

For categories of the kind $\Th \sp 0$, the condition (1) is always satisfied and the condition (2) is satisfied
if there exists such $\Th \sp 0$-algebra $H$  that either $H$ generates the given variety or the set of free
generates of $H$ has not less elements than arities of all operations.

\end{document}